\documentclass[preprint,12pt]{elsarticle}

\usepackage{amssymb}

\usepackage{amsmath}
\usepackage{multirow}
\numberwithin{equation}{section}
\newcommand\numberthis{\addtocounter{equation}{1}\tag{\theequation}}

\usepackage{url} 

\begin{document}
	
\begin{frontmatter}



\title{An integral-like numerical approach for solving Burgers' equation}


\author{Somrath Kanoksirirath\corref{cor1}\fnref{label1}}
\affiliation[label1]{organization={NSTDA Supercomputer Center (ThaiSC), National Science and Technology Development Agency (NSTDA)},
	city={Pathum Thani},
	postcode={12120}, 
	country={Thailand}}

\begin{abstract}
An unconventional approach is applied to solve the one-dimensional Burgers' equation. 
It is based on spline polynomial interpolations and Hopf-Cole transformation.
Taylor expansion is used to approximate the exponential term in the transformation, then the analytical solution of the simplified equation is discretized to form a numerical scheme, involving various special functions.
The derived scheme is explicit and adaptable for parallel computing.
However, some types of boundary condition cannot be specified straightforwardly.
Three test cases were employed to examine its accuracy, stability, and parallel scalability.
In the aspect of accuracy, the schemes employed cubic and quintic spline interpolation performs equally well, managing to reduce the $\ell_{1}$, $\ell_{2}$ and $\ell_{\infty}$ error norms down to the order of $10^{-4}$.
Due to the transformation, their stability condition $\nu \Delta t/\Delta x^2 > 0.02$ includes the viscosity/diffusion coefficient $\nu$.
From the condition, the schemes can run at a large time step size $\Delta t$ even when grid spacing $\Delta x$ is small.
These characteristics suggest that the method is more suitable for operational use than for research purposes.
\end{abstract}




\end{frontmatter}

\section{Introduction}
\label{Sec:Intro}

The partial differential equation (PDE) of the form
\begin{equation}
	\frac{\partial u}{\partial t} + u \frac{\partial u}{\partial x} = \nu \frac{\partial^2 u}{\partial x^2}
	\label{eq:burger}
\end{equation}
is called one-dimensional Burgers' equation. 
It models physical transport phenomena in fluid flows, turbulence, traffic flows and other areas \cite{reviewBurger}. 
It is a vital part of Navier-Stokes equations where $u$ represents the velocity of the fluid and $\nu$ is a positive viscosity/diffusion constant.
The equation forms the core of computational fluid dynamics, weather models, ocean models, and hydrodynamic models.

The second term, $u u_x$, is non-linear, which hinders the development of simple, stable and accurate numerical methods for solving various physics equations in practice.
This study proposes an approach to address this non-linear term while keeping the scheme simple and explicit, as highly complex schemes often suffer from poor parallel scalability, limiting their practical use.
Additionally, in contract to most explicit schemes that have severe stability conditions resulting in higher computational costs, our scheme is stable even at large time step sizes.

The solution of Burgers' equation can be obtained analytically by converting it to a diffusion equation, $\phi_t = \nu \phi_{xx}$, using Hopf-Cole transform, Eqs. (\ref{eq:HopfCole_forward}) and (\ref{eq:HopfCole_backward}), and then solving it using Fourier transform or other methods. 
However, this approach is inconvenient and unsuitable for arbitrary initial and boundary conditions commonly encountered in practice.
Therefore, several numerical schemes have been developed.
\begin{align}
	\phi(x,t) &= \mathrm{exp}\Big(\, \frac{-1}{2\nu} \int_{0'}^{x} u(\xi, t) \, \mathrm{d}\xi \Big) 
	\label{eq:HopfCole_forward} 
	\\
	u &= -2\nu \frac{\phi_x}{\phi}
	\label{eq:HopfCole_backward}
\end{align}

The standard numerical approach to tackle the non-linear term is to linearize it by assuming that $u$ in $u u_x$ is locally constant, as done in previous works such as \cite{abdul2022,dogan2004,gana2014,hon1998,huang2010,aswin2017,yang2021}.
However, this approach usually has limitations, particularly in terms of accuracy. 
Huang and Abduwali \cite{huang2010} successfully developed an explicit and unconditionally stable scheme using this assumption.
Another approach is to rewrite the non-linear term as $(u^2)_x$ and solve it accordingly as in \cite{guo2016,gupta2021,jena2023,jiang2021,mohan2015,zhang2011}.
However, these methods often involve complicated, iterative or implicit schemes that require solving a large set of linear equations, which can limit computational speed and parallel scalability of the programs.

On the other hand, some researchers, e.g., \cite{kadal2006,kan2012,kum2019,kut1999,liao2008,mukun2015,pandy2009,sakai2005,xie2008,zhao2011}, have applied numerical Hopf-Cole transformation and solved the resulting diffusion equation instead.
Advanced schemes have been employed to accurately diffuse exponential profiles generated by the Hopf-Cole transformation, while the transformation itself is generally approximated using a finite number of terms of its infinite series form or integrated numerically using Gaussian quadrature.

In terms of numerical procedures, the relatively expensive finite element method is widely employed \cite{aksan2005,cal1982,cal1981,chai2020,dogan2004,kan2012,zhang2011,zhao2011}, as well as the Galerkin approach with cubic polynomials or others basis functions \cite{abdul2022,arora2013,gana2014,ghase2018,gupta2021,hon1998,jena2023,kum2019,singh2021,tam2016,aswin2017}.
Nevertheless, several studies have used the conventional finite difference method \cite{huang2010,jiang2021,kadal2006,kut1999,liao2008,mukun2015,yang2021}, with most of them being implicit schemes.
Additionally, Gao et al. \cite{gao2013} applied a particle-based scheme to the Burgers' equation.
To speed up implementation, Kumar et al. \cite{kumar2020} developed a hybrid predictor-corrector scheme, while an artificial neural network was adopted to accelerate a prior Galerkin approach in \cite{lara2022}. 

In this paper, we introduce an integral-like approach that mimics mathematical transformations and temporal integration to advance the numerical solution in time.
In contrast to finite discretization, we employ spline interpolation to represent gridded data as a continuous function.
The general idea of this approach is explained in Section \ref{Sec:IntegralLike}.
In Section \ref{Sec:Diffusion}, the time-stepping method used in solving the diffusion equation is described as an example and as an indispensable part of the scheme for solving Burgers' equation.
In Section \ref{Sec:HopfCole}, the numerical Hopf-Cole transformation is explained and the entire scheme is composed with additional remarks on programming aspect.
In Section \ref{Sec:ExpBurger}, the results of several numerical experiments are shown, including an example of solving Burgers'-Fisher equation.
Conclusions are presented in Section \ref{Sec:Conclusion}.

\section{Integral-like approach}
\label{Sec:IntegralLike}

The proposed integral-like approach is different from conventional methods in both numerical differentiation in space and numerical forwarding scheme in time.
To solve a PDE, our data grid contains not only the necessary field variables but also their derivatives, so that the variables can be represented as a continuous function using spline polynomial interpolation between adjacent grid points.
To advance the variables in time, mathematical procedures along with additional approximations, based on the corresponding analytical solution, are emulated using the known continuous spline polynomial function.
The derivatives also need to be updated in time. 
Initializing the derivatives using a finite difference method has been found to be sufficient.
Since we can split terms in a PDE and solve more but simpler PDEs successively, as discussed in \cite{fluidgraphic}, the approach can be applied to any Burger's-type equation.
In Section \ref{Sec:ExpBurger}, a Burgers'-Fisher equation, $u_t + u u_x = \nu u_{xx} - 3 u (1-u) (1-2u)$, is split and solved as an example.

In this paper, we investigate integral-like methods that utilize linear, cubic and quintic spline interpolations. 
These methods will be referred to as linear grid scheme (LG), cubic grid scheme (CG) and quintic grid scheme (QG) respectively, as in Table \ref{table:schemes}.
For the cubic and quintic schemes, the second-order central finite difference is used to initialize $u_x$ and $u_{xx}$, except at the two end points where the first-order forward and backward finite difference are employed instead.

To forward the numerical solution of the Burgers' equation in time, the Hopf-Cole transformation is applied numerically, the resulting diffusion equation is solved in the Hopf-Cole space and then converted back to the normal space. 
The application of the integral-like approach to the diffusion equation is discussed first to familiarize readers with the general idea of the method.
It is worth noting that the integral-like approach is an explicit scheme that incorporates mathematical formulas and is equivalent to a Semi-Lagrangian method when applied to the linear advection equation.

\begingroup
\setlength{\tabcolsep}{8pt} 
\renewcommand{\arraystretch}{1.25} 
\begin{center}
\begin{table*}
\centering
\caption{Integral-like schemes with their associated spline polynomial interpolations $P_{j}$ for $u,u_x$ and $u_{xx}$ between $x_{j}$ and $x_{j+1}$ where $y = x - x_{j}$ and $0\leq y < \Delta x$}
\begin{tabular}{lll}
	\hline
	\textbf{Scheme}  & \textbf{Local spline polynomial} $\mathbf{P_{j}(y)}$ & \textbf{Data grid} \\
	\hline
	$\mathrm{Linear\,(LG)}$  
	            & $a_jy + b_j$               & $u$   \\
	\hline
	$\mathrm{Cubic\,(CG)}$   
	            & $a_jy^3+b_jy^2+c_jy+d_j$   & $u$   \\
	            & $3a_jy^2+2b_jy+c_j$        & $u_x$ \\
    \hline
	$\mathrm{Quintic\,(QG)}$ 
	            & $a_jy^5+b_jy^4+c_jy^3+d_jy^2+e_jy+f_j$ & $u$   \\
	            & $5a_jy^4+4b_jy^3+3c_jy^2+2d_jy+e_j$    & $u_x$ \\
	            & $20a_jy^3+12b_jy^2+6c_jy+2d_j$         & $u_{xx}$ \\
	\hline
\end{tabular}
\label{table:schemes}
\end{table*}
\end{center}
\endgroup

\section{Method for linear diffusion}
\label{Sec:Diffusion}
The linear diffusion equation, $\phi_t = \nu \phi_{xx}$, where $\nu$ is a positive diffusion constant, can be solved analytically using the Fourier transform. 
This approach is discussed in textbooks such as \cite{olver2013} and \cite{stone2009}.
The analytical solution is
\begin{equation}
	\phi(x,t)
	=
	\int_{-\infty}^{\infty} \frac{1}{\sqrt{4\pi \nu t}} \, \mathrm{exp}\Big(-\frac{(x-\xi)^2}{4\nu t}\Big) \, \phi(\xi,0) \ \mathrm{d}\xi
	\label{eq:diffusion}
\end{equation}
with analytical boundary conditions $\phi(\infty,t) = \phi(-\infty,t) = 0$. 

The time-stepping scheme of the integral-like approach can be derived from Eq. (\ref{eq:diffusion}) by substituting the spline polynomial function $P_{j}(y,t)$ and changing the time interval. 
In the case of the cubic scheme (CG) with equally-spacing grid points, Eq. (\ref{eq:diffusion}) becomes
\begin{align*}
	\phi(x_i,t+\Delta t)
 	&=
 	\sum_{j=-\infty}^{\infty}
	 \int_{0}^{\Delta x}
  	 \frac{1}{\sqrt{4\pi \nu \Delta t}} \, \mathrm{exp}\Big(-\frac{(x_i-(\xi_j+y))^2}{4\nu \Delta t}\Big) 
	 P_j(y,t)
	 \, \mathrm{d}(\xi_j+y) \\
 	&\approx
 	\frac{1}{\sqrt{4\pi \nu \Delta t}} \,
 	\sum_{|\ell_{i,j}| \leq 5\sigma}
	\int_{0}^{\Delta x} \mathrm{exp}\Big(-\frac{(y+\ell_{i,j})^2}{4\nu \Delta t}\Big) 
	\Big( a_jy^3+b_jy^2+c_jy+d_j \Big)
	\, \mathrm{d}y
	\numberthis \label{eq:diffusion_timescheme0}
\end{align*}
Due to the Gaussian decay term, the summation range can be limited to $j$ such that the relative distance $\ell_{i,j} \equiv \xi_j - x_i$ is within a margin of five standard deviations, i.e., $5\sigma = 5\sqrt{2\nu\Delta t}$, which is chosen for this paper.

It was found experimentally that all three methods are numerically stable as long as the maginal range of $5\sigma$ is larger than the grid spacing $\Delta x$.
This is evidenced in Figure \ref{fig:Stability_tophat}.
Mathematically, the condition is equivalent to $d \geq 0.02$, where $d=\nu (\Delta t)/(\Delta x)^2$ is a non-dimensional diffusion number.
Interestingly, the stability condition for these methods is reversed compared to the stability condition of explicit finite difference schemes.
Since the $5\sigma$ length is an estimation, the value $0.02$ is not exact.
Notably, the stability condition is independent of $u$. 

Due to the summation of $j$ within the marginal range, the time complexity of the integral-like method is $\mathcal{O}( (5\sigma/\Delta x)\, n_x n_t)$ $\sim$ $\mathcal{O}((\sqrt{\nu \Delta t}/\Delta x)\, n_x n_t)$ $\sim$ $\mathcal{O}( \sqrt{\nu} \, n_x^2 n_t^{1/2})$, where $n_x$ is the number of grid points and $n_t$ is the total number of time step.
Evidently, the complexity of the algorithm is not linear.

Another implication of the marginal range is that the boundary conditions may have to be specified by a small set of points, instead of a single point exactly at the boundary.
For example, in the case of periodic boundary condition, the required outside point $-j$ on the left of the considered domain corresponds to the inside point $n-j$.
Similarly, the values at the outside point $n_x+j$ on the right are that of the $j$-th point.
For Dirichlet and no-flux boundary conditions, reflected points or their mirror images can be used.
Adding extra grid points to both ends can keep the implementation simple and is adopted here, as it also matches with the domain decomposition method for parallelization.

Next, to quantify Eq. (\ref{eq:diffusion_timescheme0}), the indefinite integral of the form, $\int y^m\,\mathrm{exp}(-(y+\ell)^2/\delta) \, \mathrm{d}y$, where $m$ is a positive integer, is evaluated recursively using Eqs. (\ref{eq:PG_m}) - (\ref{eq:PG_1}), which derived by applying integration by parts.
\begin{align*}
	PG_m 
	&\equiv \int y^m\,\mathrm{exp}\Big(-\frac{(y+\ell)^2}{\delta}\Big) \mathrm{d}y 
	\\ 
	&=
	(m-1)\,\frac{\delta}{2} \, PG_{m-2} - \ell \, PG_{m-1} -\frac{\delta}{2} \, y^{m-1} \, \mathrm{exp}\Big(-\frac{(y-\ell)^2}{\delta}\Big) \numberthis \label{eq:PG_m}
\end{align*}
where
\begin{align*}
	PG_0
	&= \frac{\sqrt{\pi\delta}}{2} \, \Bigg[ \mathrm{erf}\bigg(\frac{y+\ell}{\sqrt{\delta}}\bigg) - \mathrm{erf}\bigg(\frac{\ell}{\sqrt{\delta}}\bigg) \Bigg] \numberthis\label{eq:PG_0}\\
	PG_1
	&= -\frac{\ell \sqrt{\pi\delta}}{2} \, \Bigg[ \mathrm{erf}\bigg(\frac{y+\ell}{\sqrt{\delta}}\bigg) - \mathrm{erf}\bigg(\frac{\ell}{\sqrt{\delta}}\bigg) \Bigg]
	- \frac{\delta}{2} \, \Bigg[ \mathrm{exp}\bigg(-\frac{(y+\ell)^2}{\delta}\bigg) - \mathrm{exp}\bigg(-\frac{\ell^2}{\delta}\bigg) \Bigg]
	\numberthis\label{eq:PG_1}
\end{align*}
On the other hand, the coefficients of the cubic spline polynomial function $P_j(y)$, i.e., $a_j, b_j, c_j$ and $d_j$, are found by solving Eqs. (\ref{eq:coeff1}) - (\ref{eq:coeff4}) at each time step.
These coefficients are analytically solved beforehand to speed up the program.
\begin{align}
	P_j(\Delta x, t)
	= \phi_{j+1} 
	&= a_j (\Delta x)^3 + b_j (\Delta x)^2 + c_j (\Delta x) + d_j \label{eq:coeff1}
	\\
	\partial_x P_j(\Delta x, t)
	= \partial_x \phi_{j+1} 
	&= 3 a_j (\Delta x)^2 + 2 b_j (\Delta x) + c_j 
	\\
	P_j(0, t)
	= \phi_{j} 
	&= d_j 
	\\
	\partial_x P_j(0, t) 
	= \partial_x \phi_{j} 
	&= c_j 
	\label{eq:coeff4}
\end{align}
We denote $\partial_x \phi_{j}$ as the first derivative of $\phi$ at grid point $j$. 
These derivatives are stored in a data grid and are updated by using Eq. (\ref{eq:diffusion_divtimescheme0}), which derived by differentiating Eq. (\ref{eq:diffusion_timescheme0}) with respect to $x_i$.
\begin{align*}
	\frac{\partial}{\partial x_i} \, \phi(x_i,t+\Delta t)
	&=
	\frac{1}{\sqrt{4\pi \nu \Delta t}} \,
	\sum_{|\ell_{i,j}| \leq 5\sigma}
	\int_{0}^{\Delta x} 
	\Bigg[\frac{\partial (\xi_j - x_i)}{\partial x_i }\frac{\partial}{\partial \ell_{i,j}} \mathrm{exp}\Big(-\frac{(y+\ell_{i,j})^2}{4\nu \Delta t}\Big) \Bigg]
	\\ &\phantom{=} 
	\Big( a_jy^3+b_jy^2+c_jy+d_j \Big)
	\, \mathrm{d}y \\
	&=
	\frac{1}{\sqrt{4\pi \nu \Delta t}} \,
	\sum_{|\ell_{i,j}| \leq 5\sigma}
	\int_{0}^{\Delta x} 
	\Bigg[(-1)\, \Big(\frac{-2 (y+\ell_{i,j})}{4\nu \Delta t}\Big) \mathrm{exp}\Big(-\frac{(y+\ell_{i,j})^2}{4\nu \Delta t}\Big) \Bigg]	
	\\ &\phantom{=} 
	\Big( a_jy^3+b_jy^2+c_jy+d_j \Big)
	\, \mathrm{d}y \\
	&=
	\frac{1}{(2\nu \Delta t) \sqrt{4\pi \nu \Delta t}} \,
	\sum_{|\ell_{i,j}| \leq 5\sigma}
	\int_{0}^{\Delta x} 
	\mathrm{exp}\Big(-\frac{(y+\ell_{i,j})^2}{4\nu \Delta t}\Big)
	\\&\phantom{=} 
	\Big( a_j y^4 + (b_j+\ell_{i,j}a_j)y^3 + (c_j+\ell_{i,j}b_j)y^2 + (d_j+\ell_{i,j}c_j)y + \ell_{i,j}d_j\Big)
	\, \mathrm{d}y
	\numberthis \label{eq:diffusion_divtimescheme0}
\end{align*}
Hence, Eqs. (\ref{eq:diffusion_timescheme0}) and (\ref{eq:diffusion_divtimescheme0}) together form the complete time-stepping method for the integral-like cubic scheme.
Their computation is facilitated by Eqs. (\ref{eq:PG_m}) - (\ref{eq:coeff4}).
The derivation of the linear scheme and the quintic scheme follows a similar approach as the cubic scheme, but with a different form of spline polynomial $P_j(y)$ substituted in, and with a different number of data grids to be stored and iterated in time.

\section{Method for Hopf-Cole transformation}
\label{Sec:HopfCole}
For the cubic scheme, the first step in performing the Hopf-Cole transformation, changing $u$ to $\phi$ (Eq. \ref{eq:HopfCole_forward}), is to compute the integral of $u$ in Eq. (\ref{eq:HopfCole_forward_integrate}).
\begin{align*}
	\mathrm{Int}(u)
	\equiv 
	\int_{0'}^{x_i} u(\xi, t) \, \mathrm{d}\xi
	&= 
 	\sum_{j=0}^{i-1} \int_{0}^{\Delta x} 
	P_j(\xi, t) \, \mathrm{d}\xi 
	\\
	&= 
	\sum_{j=0}^{i-1} \int_{0}^{\Delta x} 
	\Big(a_j \xi^3 + b_j \xi^2 + c_j \xi + d_j \Big) \, \mathrm{d}\xi  
	\\
	&=
	\sum_{j=0}^{i-1}
	\ a_j \frac{(\Delta x)^4}{4} + b_j \frac{(\Delta x)^3}{3}  + c_j \frac{(\Delta x)^2}{2} + d_j (\Delta x)
	\numberthis \label{eq:HopfCole_forward_integrate}
\end{align*}
Then, the exponent $-(1/2\nu)
 \, \mathrm{Int}(u)$ in Eq. (\ref{eq:HopfCole_forward}) is represented as a quintic spline polynomial function using $\mathrm{Int}(u)$, $u$ and $u_x$.
The quintic coefficients are found by solving six linear equations developed analogously to Eqs. (\ref{eq:coeff1}) - (\ref{eq:coeff4}).

Then, for a quintic polynomial $Q_i$ of the $i$-th segment, the value of $\phi$ between $x_i$ and $x_{i+1}$ is written as a Taylor series $T_i$ of $\mathrm{exp}(Q_i)$. 
In our implementation, the number of terms in the Taylor series are varied adaptively, up to $16$ terms, to ensure that the deviation of $T_i$ is less than $0.01\%$.
The resulting $\phi$ or $T_i$, which satisfying the diffusion equation, is then solved using the method discussed in Section \ref{Sec:Diffusion}.

On the other hand, the Hopf-Cole transformation from $\phi$ to $u$ in Eq. (\ref{eq:HopfCole_backward}) is done by direct substitution.
For example, $u$ and $u_x$ of the cubic scheme are deduced from $\phi$, $\phi_x$ and $\phi_{xx}$ by using Eqs. (\ref{eq:HopfCole_backward_numer}).
\begin{equation}
	u_i = -2\nu \, \bigg( \frac{\partial_x \phi_i}{\phi_i} \bigg)
	\, ,\qquad \quad
	\partial_x u_i = -2\nu \, \bigg( \frac{\phi_i (\partial^2_x \phi_i)  - (\partial_x \phi_i)^2}{\phi^2_i} \bigg)	
	\label{eq:HopfCole_backward_numer}
\end{equation}

To sum up, an integral-like scheme for solving Burgers' equation consists of the following components:
(1) Numerical Hopf-Cole integration scheme for transforming from $u$ space to $\phi$ space,
(2) Numerical diffusion scheme for advancing $\phi(x,t)$ and its required derivatives to $\phi(x,t+\Delta t)$,
(3) Numerical Hopf-Cole differentiation scheme for transforming $\phi(x,t+\Delta t)$ and its derivatives back to $u$ space, 
and (4) Spline interpolation method for representing discrete data points as a continuous function for the computation in (1)-(3).
One complete time iteration step $\Delta t$ for solving Burgers' equation is composed of (1), (2) and (3).
In addition, a finite difference scheme is employed for computing the initial derivatives such as $u_x(x,t_0)$ from $u(x,t_0)$.
Our source code is publicly available on GitHub (\url{https://github.com/SKanoksi}).

According to our implementation, because the numerical Hopf-Cole integration step involves rapid exponential decay/growth, numerical round-off error must be carefully minimized. 
At least double-precision floating-point data type may be a minimum. 
Moreover, additional memory allocation is required for the derivatives of the cubic scheme and the quintic scheme, resulting in space complexity of $\mathcal{O}(4n_x)$ and $\mathcal{O}(6n_x)$ respectively.
The common factor of $2$ is because a temporary array is needed for keeping updated value at each position, before entirely transferred to the primary array in the final step.
For our current implementation, another $n_x$ array is used to temporarily store $\mathrm{Int}(u)$.

An implementation of an adaptive grid was explored, where grid points are redistributed based on the path length of $u$, numerically approximated using Gaussian quadrature.
However, the numerical solution was found to be less accurate due to numerical errors introduced when repeatedly rearranging the grids.
Therefore, our unsuccessful implementation of the adaptive grid is briefly noted in this paper to inform the possibility that adaptive grid may not improve integral-like schemes in general.

\section{Numerical experiment}
\label{Sec:ExpBurger}
In this section, four example cases are used to evaluate the integral-like methods in Table \ref{table:schemes}.
The accuracy of the methods is tested using Example 1 and 2, while numerical stability and parallel scalability are evaluated using Example 3. 
Example 4 investigates the viability of the split approach for more complicated problems.

To quantify the accuracy of the schemes when comparing with exact/analytical solutions $f(x)$, the $\ell_1$, $\ell_2$ and $\ell_{\infty}$ error norms are adopted.
These error norms are calculated as shown below, for equally-spacing grid.
\begin{align}
	\ell_{1}(u) &= \dfrac{\sum_{j=0}^{n} |u_j-f(x_j)| }{\sum_{j=0}^{n} |f(x_j)|} 
	\, , \quad 
	\\
	\ell_{2}(u) &= \dfrac{ \sqrt{ \sum_{j=0}^{n} (u_j-f(x_j))^2 }}{\sqrt{\sum_{j=0}^{n} f(x_j)^2}} 
	\, , \quad
	\\
	\ell_{\infty}(u) &= \dfrac{ \mathrm{max} |u_j-f(x_j)|}{\mathrm{max}|f(x_j)|} 
\end{align}

\textbf{Example 1.} Burgers' equation (Eq. \ref{eq:burger}) with the initial condition
\begin{equation}
	u(x,0) = \mathrm{sin}(\pi x)
	\label{eq:ex2_init}
\end{equation}
and the exact solution 
\begin{equation}
	u(x,t) = \dfrac{4\pi\nu \ \sum_{k=1}^{\infty} \, k \, A_k \, \mathrm{sin}(k\pi x) \, \mathrm{exp}(-k^2\pi^2 \nu t) }
	{A_0 + 2\, \sum_{k=1}^{\infty} \, A_k \, \mathrm{cos}(k\pi x) \, \mathrm{exp}(-k^2\pi^2 \nu t) }
	\label{eq:ex2_ana}
\end{equation}
where
\begin{equation}
	A_k = \int_{0}^{1} \, \mathrm{cos}(k \pi x) \, \mathrm{exp}\Big( \frac{\mathrm{cos}(\pi x)-1}{2\pi \nu} \Big) \, \mathrm{d}x
\end{equation}
are considered.
Unlike in \cite{gao2013} and \cite{Sarbo2014}, periodic boundary conditions are employed.

\begin{figure}[h]
	\includegraphics[width=\linewidth]{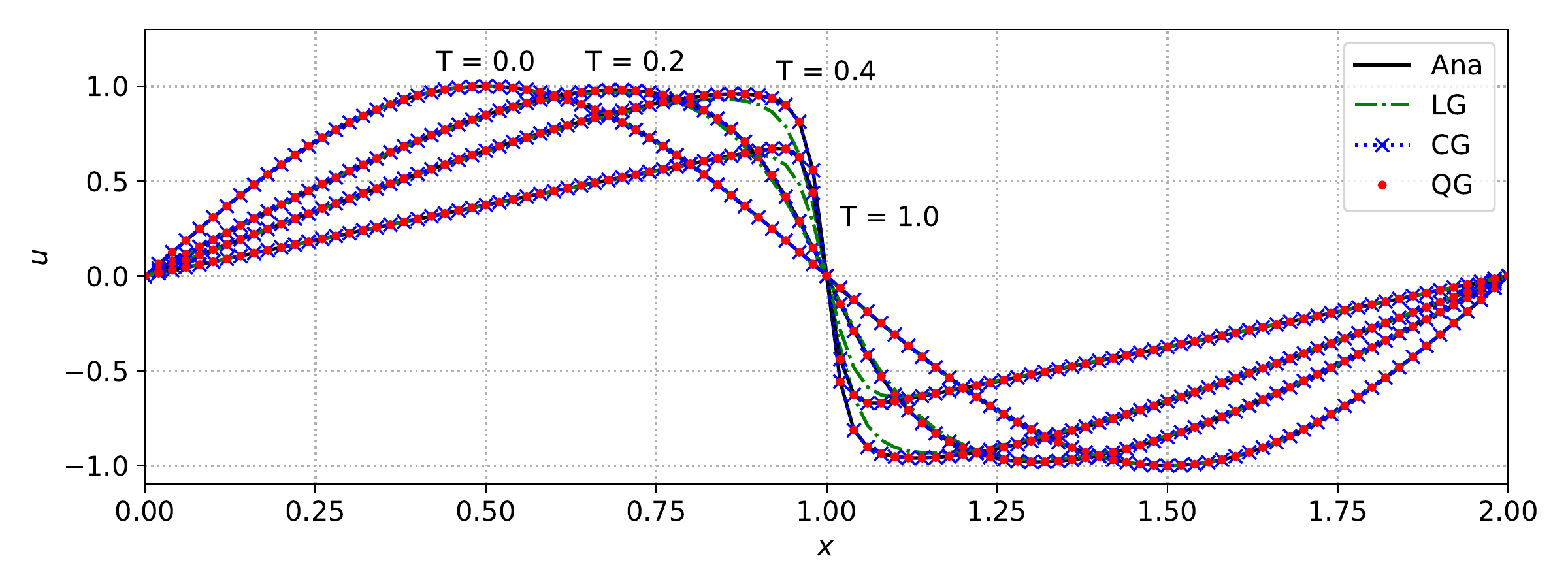}
	\caption{Numerical solutions and exact solution of Example 1: $\nu = 0.01$, $\Delta x = 0.02$, and $\Delta t = 0.005$}
	\label{fig:Sol_sine}
\end{figure}

The exact solution and numerical results are displayed in Figure \ref{fig:Sol_sine} for $\nu = 0.01$, $\Delta x = 0.02$, and $\Delta t = 0.005$, while the numerical values at some grid points are given in Table \ref{table:pos_error_sine}.
The error norms for four different pairs of $\Delta x$ and $\Delta t$ are shown in Table \ref{table:errornorm_sine}.
  
\begin{figure}[h]
	\includegraphics[width=\linewidth]{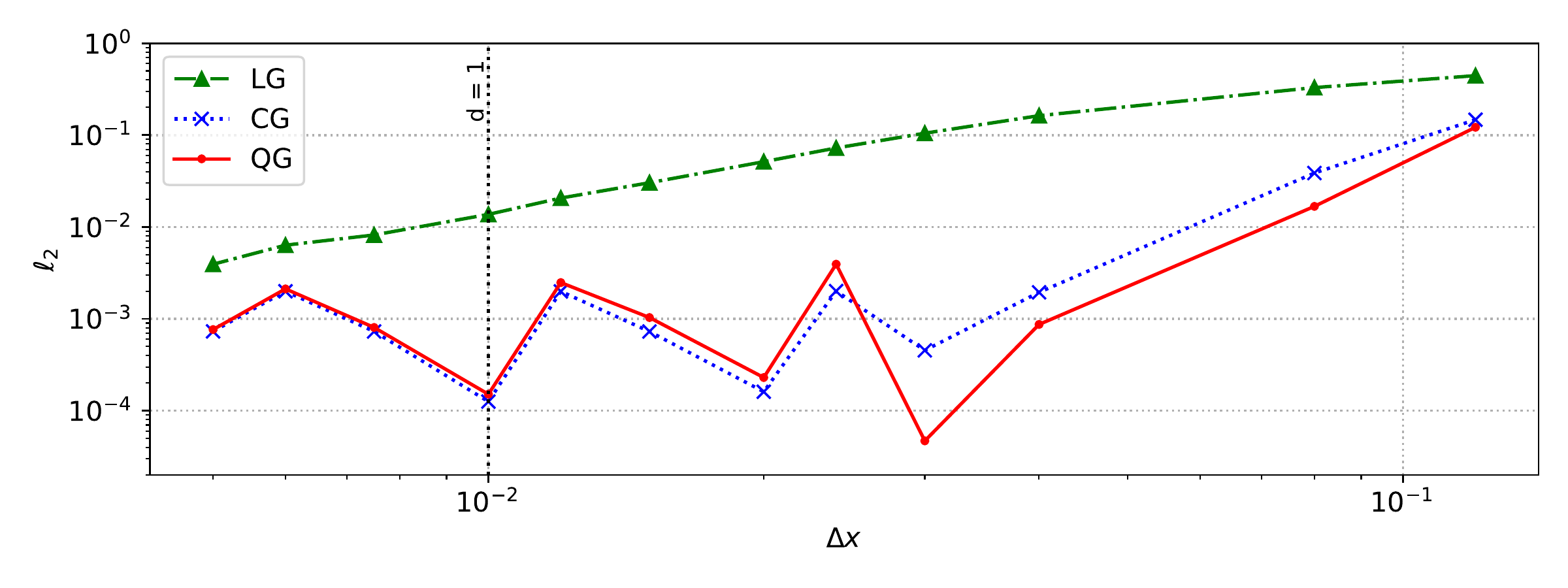}
	\caption{$\ell_2$ error norm of Example 2 as a function of $\Delta x$ using $\nu = 0.005$, $\Delta x = 0.02$, $\Delta t = 0.02$ and $n_t = 100$}
	\label{fig:Acc_pulse}
\end{figure}

From the results, the linear scheme (LG) always performs the worst, while the cubic scheme (CG) and the quintic scheme (QG) are relatively comparable.
For CG and QG, as the grid spacing and time step size are reduced, the error norms may not further decrease as seen in Figure \ref{fig:Acc_pulse}, where $\Delta x$ is varied from $0.005$ to $0.04$.
This behavior may be caused by truncation error introduced when approximating $\mathrm{exp}$ and $\mathrm{erf}$ functions, which are frequently used in the schemes.
As a result, the order of accuracy of CG and QG remains nearly unchanged, while LG is approximately second-order accurate, but probably becomes constant at around $\ell_2 = 10^{-3}$ as well.

\begin{center}
\begin{table*}
	\centering
	\caption{Numerical results and exact solution of Example 1 at $t=1.0$ using $\nu = 0.1$ and $\Delta t = 0.01$}
	\begin{tabular}{cccccc}
		\hline
		\textbf{Position} & $\Delta \mathbf{x}$ & \textbf{LG} & \textbf{CG} & \textbf{QG} & \textbf{Exact} \\
		\hline
		\multirow{2}{*}{$0.1$} & $0.02$ & $0.065558$ & $0.066306$ & $0.066305$ & \multirow{2}{*}{$0.066316$} \\
		& $0.01$ & $0.066101$ & $0.066284$ & $0.066283$ & \\
		\hline
		\multirow{2}{*}{$0.2$} & $0.02$ & $0.129578$ & $0.131189$ & $0.131187$ & \multirow{2}{*}{$0.131209$} \\
		& $0.01$ & $0.130750$ & $0.131145$ & $0.131144$ & \\
		\hline
		\multirow{2}{*}{$0.3$} & $0.02$ & $0.190039$ & $0.192755$ & $0.192752$ & \multirow{2}{*}{$0.192786$} \\
		& $0.01$ & $0.192020$ & $0.192689$ & $0.192687$ & \\
		\hline
		\multirow{2}{*}{$0.4$} & $0.02$ & $0.243809$ & $0.247999$ & $0.247995$ & \multirow{2}{*}{$0.248041$} \\
		& $0.01$ & $0.246874$ & $0.247910$ & $0.247908$ & \\
		\hline
		\multirow{2}{*}{$0.5$} & $0.02$ & $0.285749$ & $0.291864$ & $0.291859$ & \multirow{2}{*}{$0.291916$} \\
		& $0.01$ & $0.290232$ & $0.291752$ & $0.291749$ & \\
		\hline
		\multirow{2}{*}{$0.6$} & $0.02$ & $0.307656$ & $0.316005$ & $0.316000$ & \multirow{2}{*}{$0.316068$} \\
		& $0.01$ & $0.313784$ & $0.315876$ & $0.315872$ & \\
		\hline
		\multirow{2}{*}{$0.7$} & $0.02$ & $0.297804$ & $0.308022$ & $0.308017$ & \multirow{2}{*}{$0.308089$} \\
		& $0.01$ & $0.305304$ & $0.307884$ & $0.307880$ & \\
		\hline
		\multirow{2}{*}{$0.8$} & $0.02$ & $0.243417$ & $0.253657$ & $0.253653$ & \multirow{2}{*}{$0.253718$} \\
		& $0.01$ & $0.250927$ & $0.253533$ & $0.253529$ & \\
		\hline
		\multirow{2}{*}{$0.9$} & $0.02$ & $0.139279$ & $0.146027$ & $0.146025$ & \multirow{2}{*}{$0.146065$} \\
		& $0.01$ & $0.144223$ & $0.145951$ & $0.145949$ & \\
		\hline
	\end{tabular}
	\label{table:pos_error_sine}
\end{table*}
\end{center}

\begin{center}
\begin{table*}
	\centering
	\caption{Error norms of Example 1 at $t=1.0$ using $\nu = 0.1$.}
	\begin{tabular}{ccccccc}
		\hline
		\textbf{Scheme} & $\Delta\mathbf{x}$ & $\Delta\mathbf{t}$ & $\mathbf{d}$ & $\mathbf{\ell_\mathrm{1}}$ & $\mathbf{\ell_\mathrm{2}}$ & $\mathbf{\ell_{\infty}}$ \\
		\hline
		\multirow{4}{*}{$\mathrm{LG}$} & $0.02$ & $0.0200$ &  $5.0$ & $1.35\mathrm{e}-02$ & $1.46\mathrm{e}-02$ & $1.72\mathrm{e}-02$ \\
		& $0.02$ & $0.0100$ &  $2.5$ & $2.64\mathrm{e}-02$ & $2.85\mathrm{e}-02$ & $3.35\mathrm{e}-02$ \\
		& $0.01$ & $0.0100$ & $10.0$ & $7.19\mathrm{e}-03$ & $7.75\mathrm{e}-03$ & $9.07\mathrm{e}-03$ \\
		& $0.01$ & $0.0025$ &  $2.5$ & $2.70\mathrm{e}-02$ & $2.92\mathrm{e}-02$ & $3.42\mathrm{e}-02$ \\
		\hline
		\multirow{4}{*}{$\mathrm{CG}$} & $0.02$ & $0.0200$ &  $5.0$ & $4.25\mathrm{e}-04$ & $4.33\mathrm{e}-04$ & $4.56\mathrm{e}-04$ \\
		& $0.02$ & $0.0100$ &  $2.5$ & $1.97\mathrm{e}-04$ & $2.00\mathrm{e}-04$ & $2.11\mathrm{e}-04$ \\
		& $0.01$ & $0.0100$ & $10.0$ & $6.07\mathrm{e}-04$ & $6.16\mathrm{e}-04$ & $6.46\mathrm{e}-04$ \\
		& $0.01$ & $0.0025$ &  $2.5$ & $8.27\mathrm{e}-04$ & $8.39\mathrm{e}-04$ & $8.76\mathrm{e}-04$ \\
		\hline
		\multirow{4}{*}{$\mathrm{QG}$} & $0.02$ & $0.0200$ &  $5.0$ & $4.40\mathrm{e}-04$ & $4.48\mathrm{e}-04$ & $4.72\mathrm{e}-04$ \\
		& $0.02$ & $0.0100$ &  $2.5$ & $2.12\mathrm{e}-04$ & $2.15\mathrm{e}-04$ & $2.26\mathrm{e}-04$ \\
		& $0.01$ & $0.0100$ & $10.0$ & $6.18\mathrm{e}-04$ & $6.28\mathrm{e}-04$ & $6.58\mathrm{e}-04$ \\
		& $0.01$ & $0.0025$ &  $2.5$ & $8.97\mathrm{e}-04$ & $9.10\mathrm{e}-04$ & $9.50\mathrm{e}-04$ \\
		\hline
	\end{tabular}
	\label{table:errornorm_sine}
\end{table*}
\end{center}

\begin{figure}[h]
	\includegraphics[width=\linewidth]{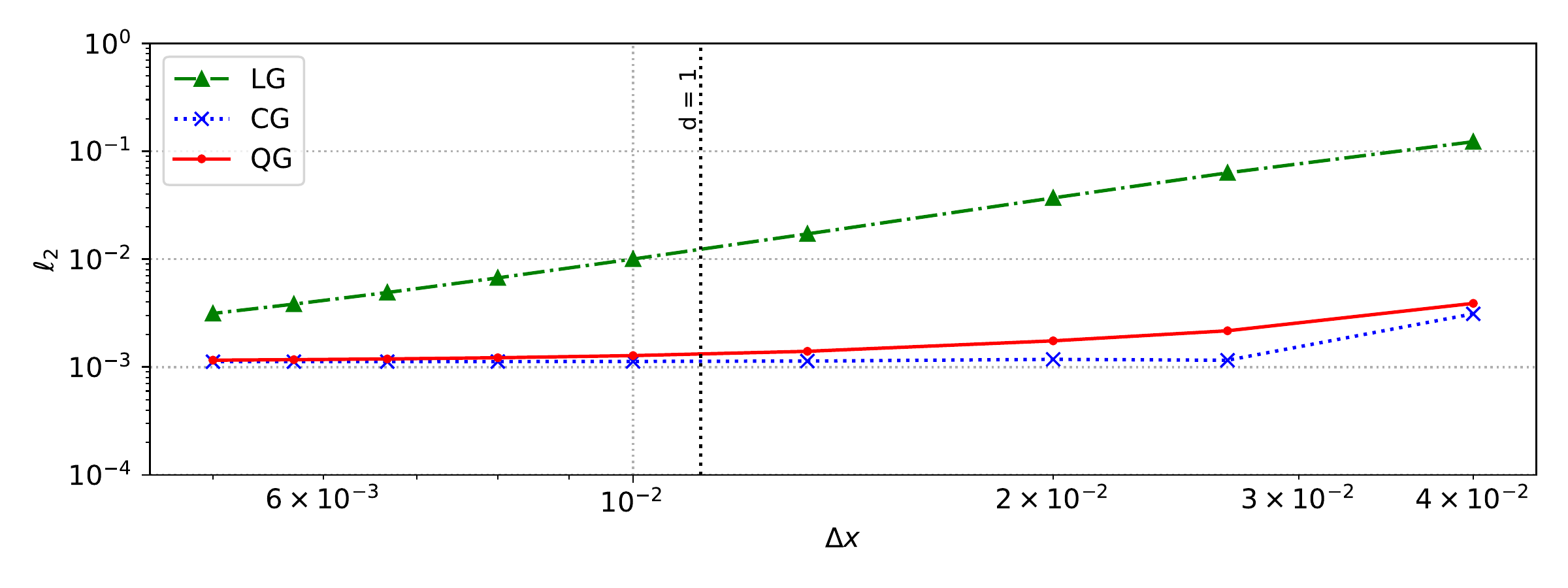}
	\caption{$\ell_2$ error norm of Example 1 as a function of $\Delta x$ using $\nu = 0.1$, $\Delta t = 0.0125$ and $n_{t} = 100$.}
	\label{fig:Acc_sine}
\end{figure}

\newpage

\textbf{Example 2.} Burgers' equation (Eq. \ref{eq:burger}) with the Dirichlet boundary condition
\begin{equation}
	u(0,t) = u(1.2,t) = 0
	\label{eq:ex1_bc}
\end{equation}
and the exact solution 
\begin{equation}
	u(x,t) = \dfrac{x/t}{1 + \sqrt{t/t_0}\ \mathrm{exp}(x^2/4\nu t)}
	\label{eq:ex1_init}	
\end{equation}
where $t \geq 1$ and $t_0 = \mathrm{exp}(0.125/\nu)$ are considered as in \cite{benton1972}, \cite{gao2013} and \cite{Sarbo2014}.
In our case, the Dirichlet boundary conditions are modeled using the inverted reflection of the solution.

\begin{figure}[h]
	\includegraphics[width=\linewidth]{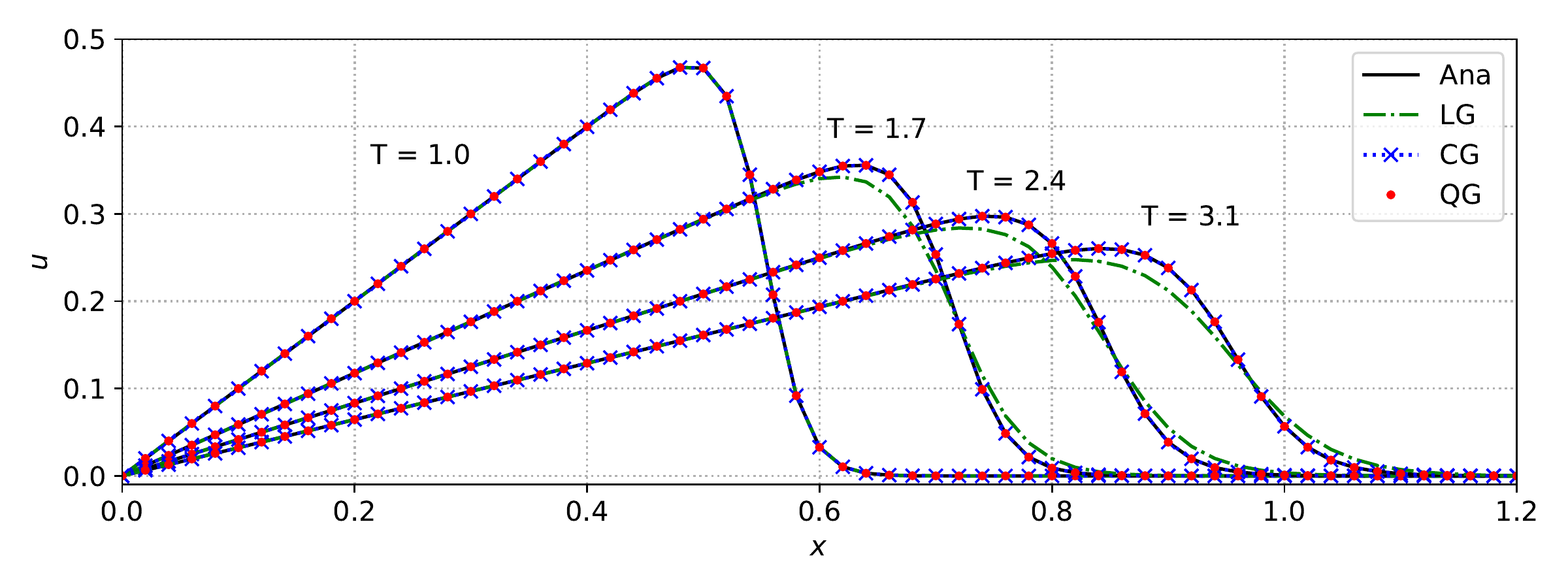}
	\caption{Numerical solutions and exact solution of Example 2: $\nu = 0.005$, $\Delta x = 0.02$ and $\Delta t = 0.02$}
	\label{fig:Sol_pulse}
\end{figure}

The results of this example, using $\nu = 0.005$, $\Delta x = 0.02$ and $\Delta t = 0.02$, are shown in Figure \ref{fig:Sol_pulse} and Table \ref{table:pos_error_pulse}, while Table \ref{table:errornorm_pulse} and Figure \ref{fig:Acc_pulse} provide more insight on their order of accuracy.
The results in Figure \ref{fig:Acc_pulse} exhibit similar features as in Figure \ref{fig:Acc_sine} of previous example, except for the oscillation of CG and QG at small $\Delta x$.
The swing may be due to the accumulation of errors, similar to Runge's and/or Gibbs phenomenon near the steep slope.
This behavior does not occur in Example 1, probably because the abrupt change is stationary at the grid point $x=0$, while, in this example, the change is not always on a grid points as $\Delta x$ is varied.

\begin{center}
\begin{table*}
	\centering
	\caption{Numerical results and exact solution of Example 2 at $t=2.4$ using $\nu = 0.005$ and $\Delta t = 0.02$}
	\begin{tabular}{cccccc}
		\hline
		\textbf{Position} & $\Delta\mathbf{x}$ & \textbf{LG} & \textbf{CG} & \textbf{QG} & \textbf{Exact} \\
		\hline
		\multirow{2}{*}{$0.1$} & $0.02$ & $0.041664$ & $0.041663$ & $0.041662$ & \multirow{2}{*}{$0.041667$} \\
		                     & $0.01$ & $0.041665$ & $0.041665$ & $0.041665$ & \\
        \hline
		\multirow{2}{*}{$0.2$} & $0.02$ & $0.083328$ & $0.083328$ & $0.083326$ & \multirow{2}{*}{$0.083333$} \\
		                     & $0.01$ & $0.083331$ & $0.083331$ & $0.083330$ & \\
        \hline
		\multirow{2}{*}{$0.3$} & $0.02$ & $0.124995$ & $0.124997$ & $0.124994$ & \multirow{2}{*}{$0.125000$} \\
		                     & $0.01$ & $0.124996$ & $0.124996$ & $0.124996$ & \\
        \hline
		\multirow{2}{*}{$0.4$} & $0.02$ & $0.166648$ & $0.166661$ & $0.166656$ & \multirow{2}{*}{$0.166665$} \\
		                     & $0.01$ & $0.166659$ & $0.166660$ & $0.166659$ & \\
        \hline
		\multirow{2}{*}{$0.5$} & $0.02$ & $0.208207$ & $0.208312$ & $0.208305$ & \multirow{2}{*}{$0.208318$} \\
		                     & $0.01$ & $0.208299$ & $0.208311$ & $0.208309$ & \\
        \hline
		\multirow{2}{*}{$0.6$} & $0.02$ & $0.248953$ & $0.249808$ & $0.249797$ & \multirow{2}{*}{$0.249816$} \\
		                     & $0.01$ & $0.249681$ & $0.249806$ & $0.249804$ & \\
        \hline
		\multirow{2}{*}{$0.7$} & $0.02$ & $0.281531$ & $0.288475$ & $0.288445$ & \multirow{2}{*}{$0.288472$} \\
		                     & $0.01$ & $0.287090$ & $0.288458$ & $0.288454$ & \\
        \hline
		\multirow{2}{*}{$0.8$} & $0.02$ & $0.239172$ & $0.266132$ & $0.266155$ & \multirow{2}{*}{$0.266228$} \\
		                     & $0.01$ & $0.258598$ & $0.266184$ & $0.266179$ & \\
        \hline
		\multirow{2}{*}{$0.9$} & $0.02$ & $0.055092$ & $0.038663$ & $0.038622$ & \multirow{2}{*}{$0.038651$} \\
		                     & $0.01$ & $0.043197$ & $0.038638$ & $0.038633$ & \\
        \hline
		\multirow{2}{*}{$1.0$} & $0.02$ & $0.003517$ & $0.000902$ & $0.000911$ & \multirow{2}{*}{$0.000912$} \\
							 & $0.01$ & $0.001385$ & $0.000911$ & $0.000912$ & \\
        \hline
		\multirow{2}{*}{$1.1$} & $0.02$ & $0.000137$ & $0.000012$ & $0.000013$ & \multirow{2}{*}{$0.000013$} \\
							 & $0.01$ & $0.000026$ & $0.000013$ & $0.000013$ & \\
		\hline
	\end{tabular}
	\label{table:pos_error_pulse}
\end{table*}
\end{center}

\begin{center}
\begin{table*}
	\centering
	\caption{Error norms of Example 2 at $t=2.4$, using $\nu = 0.005$.}
	\begin{tabular}{ccccccc}
		\hline
		\textbf{Scheme} & $\Delta\mathbf{x}$ & $\Delta\mathbf{t}$ & $\mathbf{d}$ & $\mathbf{\ell_\mathrm{1}}$ & $\mathbf{\ell_\mathrm{2}}$ & $\mathbf{\ell_{\infty}}$ \\
		\hline
		\multirow{4}{*}{$\mathrm{LG}$} & $0.02$ & $0.035$ & $0.4375$ & $1.77\mathrm{e}-02$ & $2.86\mathrm{e}-02$ & $5.55\mathrm{e}-02$ \\
		                    & $0.02$ & $0.020$ & $0.2500$ & $3.09\mathrm{e}-02$ & $4.85\mathrm{e}-02$ & $9.10\mathrm{e}-02$ \\
		                    & $0.01$ & $0.020$ & $1.0000$ & $7.82\mathrm{e}-03$ & $1.30\mathrm{e}-02$ & $2.56\mathrm{e}-02$ \\
		                    & $0.01$ & $0.005$ & $0.2500$ & $3.12\mathrm{e}-02$ & $4.91\mathrm{e}-02$ & $9.18\mathrm{e}-02$ \\
        \hline
        \multirow{4}{*}{$\mathrm{CG}$} & $0.02$ & $0.035$ & $0.4375$ & $2.78\mathrm{e}-04$ & $3.95\mathrm{e}-04$ & $8.24\mathrm{e}-04$ \\
		                    & $0.02$ & $0.020$ & $0.2500$ & $8.78\mathrm{e}-05$ & $1.48\mathrm{e}-04$ & $3.98\mathrm{e}-04$ \\
		                    & $0.01$ & $0.020$ & $1.0000$ & $6.89\mathrm{e}-05$ & $9.10\mathrm{e}-05$ & $1.77\mathrm{e}-04$ \\
		                    & $0.01$ & $0.005$ & $0.2500$ & $1.71\mathrm{e}-04$ & $2.01\mathrm{e}-04$ & $3.86\mathrm{e}-04$ \\
        \hline
        \multirow{4}{*}{$\mathrm{QG}$} & $0.02$ & $0.035$ & $0.4375$ & $3.91\mathrm{e}-04$ & $5.41\mathrm{e}-04$ & $1.01\mathrm{e}-03$ \\
		                    & $0.02$ & $0.020$ & $0.2500$ & $1.31\mathrm{e}-04$ & $1.67\mathrm{e}-04$ & $3.04\mathrm{e}-04$ \\
		                    & $0.01$ & $0.020$ & $1.0000$ & $8.39\mathrm{e}-05$ & $1.09\mathrm{e}-04$ & $2.01\mathrm{e}-04$ \\
		                    & $0.01$ & $0.005$ & $0.2500$ & $3.31\mathrm{e}-04$ & $3.81\mathrm{e}-04$ & $6.41\mathrm{e}-04$ \\
		\hline
	\end{tabular}
	\label{table:errornorm_pulse}
\end{table*}
\end{center}

\newpage

\begin{figure}[h]
	\includegraphics[width=\linewidth]{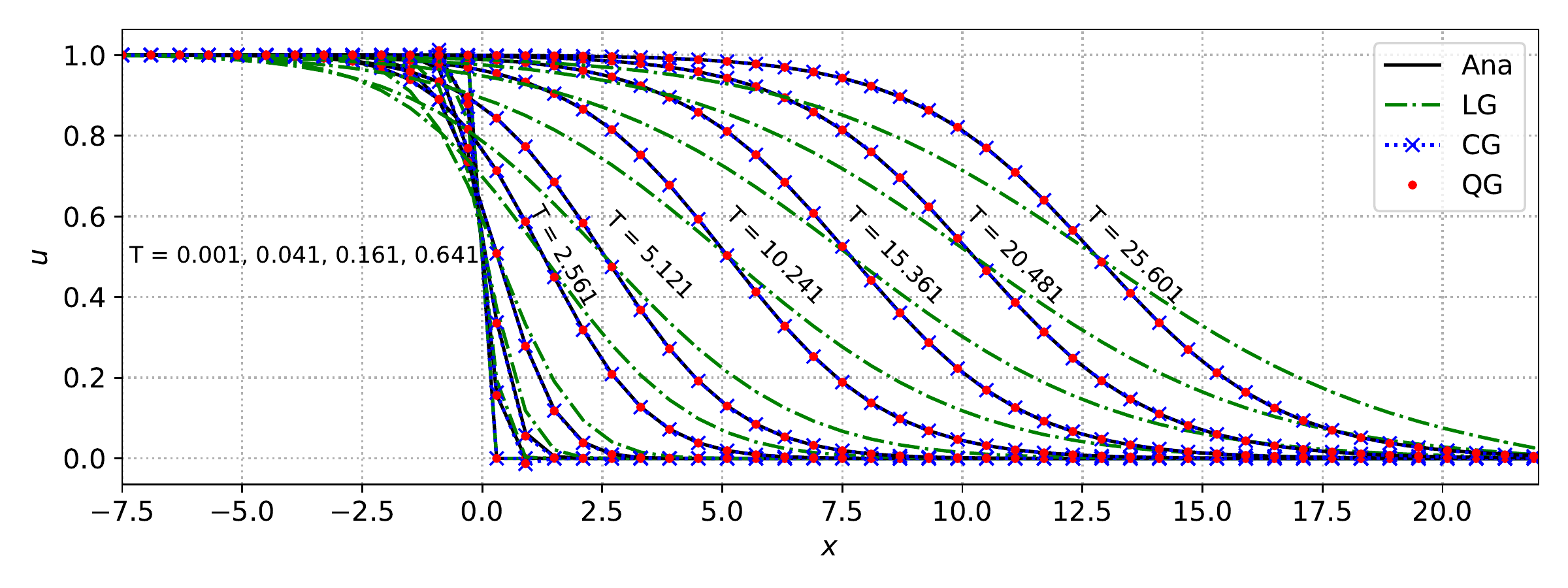}
	\caption{Numerical solutions and exact solution of Example 3: $\nu = 1$, $\Delta x = 0.6$ and $\Delta t = 0.04$}
	\label{fig:Sol_tophat_travel1}
\end{figure}

\begin{figure}[h]
	\includegraphics[width=\linewidth]{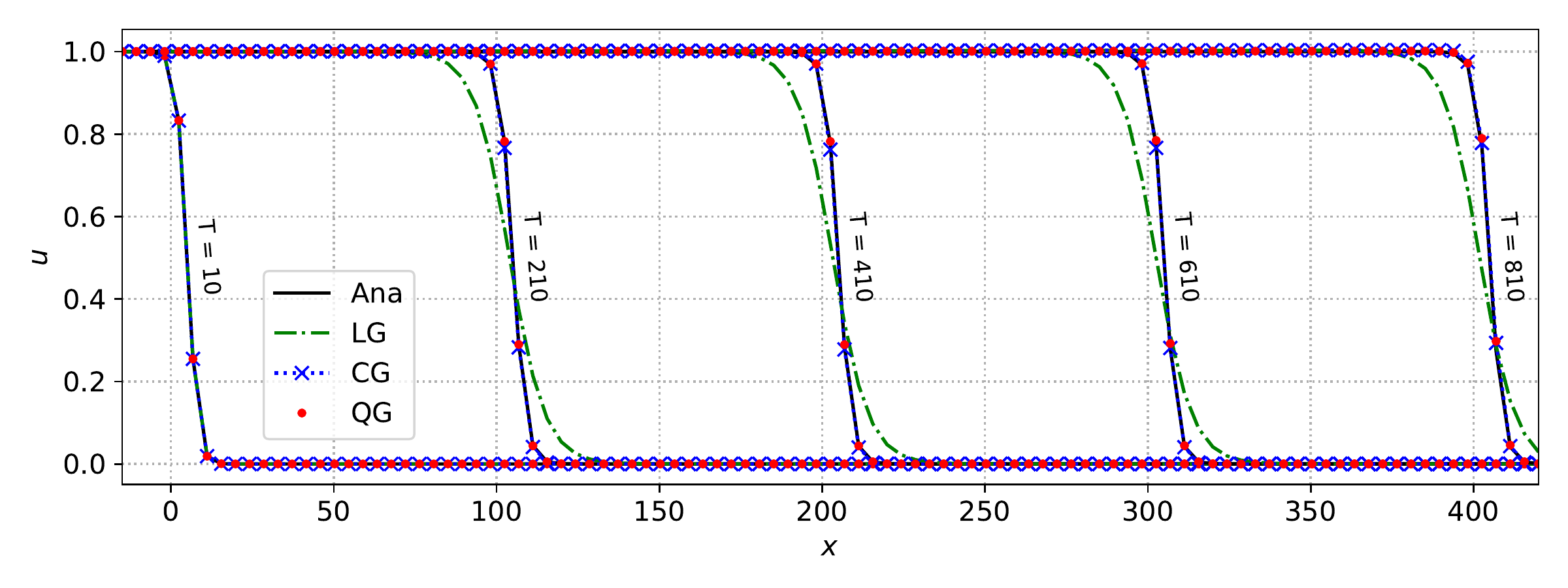}
	\caption{Numerical solutions and exact solution of Example 3: $\nu = 1$, $\Delta x = 4.35$ and $\Delta t = 0.8$}
	\label{fig:Sol_tophat_travel3}
\end{figure}

\textbf{Example 3.} Burgers' equation (Eq. \ref{eq:burger}) with a step initial condition at $x=0$ is used to study both the stability and parallel scalability of the integral-like method. 
The exact solution for infinite boundary conditions $u(-\infty, t)=1$ and $u(\infty, t)=0$ is
\begin{equation}
	u(x,t) = \dfrac{\mathrm{erfc}\big(\frac{x-t}{2\sqrt{\nu t}}\big)}
	{\mathrm{erfc}\big(-\frac{x}{2\sqrt{\nu t}}\big) \ \mathrm{exp}\big( \frac{x-t/2}{2\nu}\big) + \mathrm{erfc}\big(\frac{x-t}{2\sqrt{\nu t}}\big)}
	\label{eq:ex3_ana}
\end{equation}
which becomes a traveling-wave solution thereafter about $t=10$.
Using $\nu = 1$, $\Delta x = 0.6$ and $\Delta t = 0.04$, the results are shown in Figure \ref{fig:Sol_tophat_travel1} and \ref{fig:Sol_tophat_travel3} for different experimental setups. 

\begin{figure}[h]
	\includegraphics[width=\linewidth]{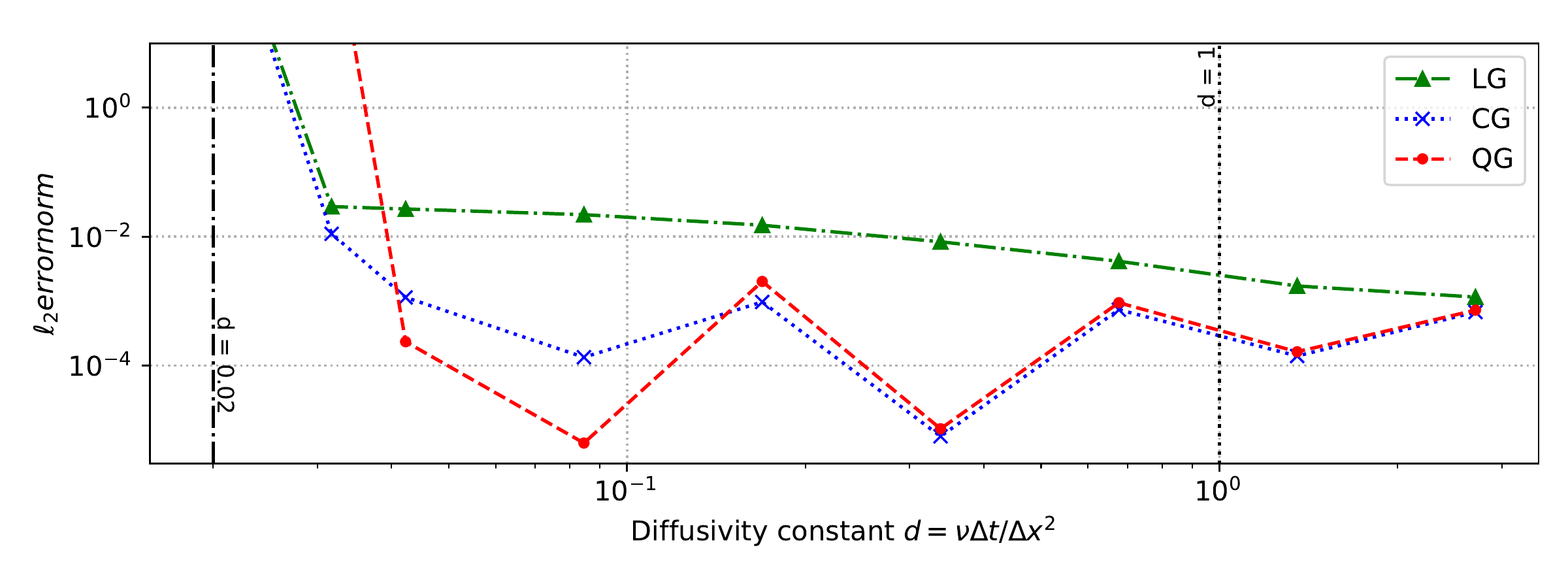}
	\caption{$\ell_2$ error norm of Example 3 as a function of $\nu$ using $\Delta x = 4.35$, $n_x = 401$, $\Delta t = 0.8$ and $n_t = 100$}
	\label{fig:Stability_tophat}
\end{figure}

In Figure \ref{fig:Stability_tophat}, the $\ell_2$ error norm is plotted against the non-dimensional diffusion number $d = \nu \Delta t/\Delta x^2$.
It can be seen that the stability condition of the integral-like method is not $d<1$, as also indicated by Figure \ref{fig:Acc_sine} and \ref{fig:Acc_pulse}, but rather $d>0.02$.
This finding demonstrates that, unlike most simple explicit methods, the integral-like method is stable at large $\Delta t$.
On the other hand, the condition $d>0.02$ implies that the marginal range $5\sigma$ has to be larger than the grid spacing, as discussed in Section \ref{Sec:Diffusion}.
The remaining problem encountered when $\nu$ is very small can be resolved by simply increasing $\Delta t$.

A weak-scaling experiment was performed on the TARA cluster of the NSTDA supercomputer center (ThaiSC) to test the parallel efficiency of the integral-like method.
Both the number of grid points and the number of time steps are adjusted to scale computation workload while having the grid spacing $\Delta x$ and the marginal range $5\sigma = 5\sqrt{2\nu\Delta t}$ approximately unaltered.
In other words, roughly the same number of $u(x_j,t)$ is used in updating $u(x_i,t+\Delta t)$.

With $N_\mathrm{s}$ representing a running variable, the number of employed CPU cores is $N_\mathrm{core}=N_\mathrm{s}^2$ operating on the total number of grid points $n_x = 500\,N_\mathrm{s}+1$ for domain $[-15, 400\,N_\mathrm{s}+20]$ to run the simulation from $t_0 = 10$ to $t = 800\,N_\mathrm{s}+10$.
The time step size $\Delta t$ is varied to explore the influence of the marginal range on parallel scalability.
This is because the bigger the marginal range is, the larger the overlapped areas of the domain decomposition algorithm are.
The outputs are shown in Figure \ref{fig:Sol_tophat_travel3}.

\begin{figure}[h]
	\includegraphics[width=\linewidth]{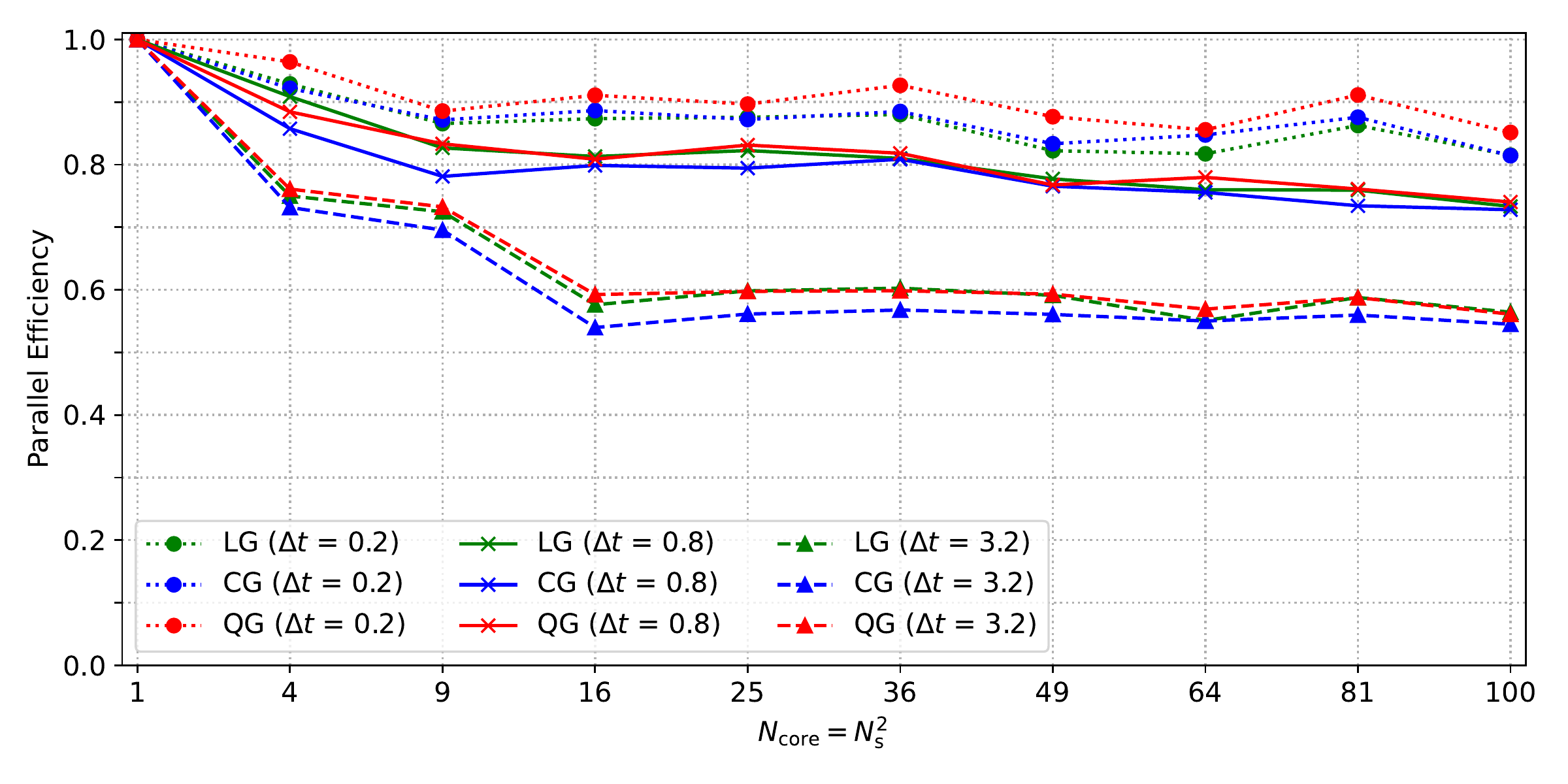}
	\caption{Weak-scaling parallel efficiency of the integral-like schemes for Burgers' equation using time step size $\Delta t = 0.2$, $0.8$ and $3.2$.}
	\label{fig:WeakSca}
\end{figure}

TARA compute nodes, equipped with two Intel Xeon Gold 6148 CPU (2.40GHz) and 192GB of RAM, are employed.
Our program was coded in Python and parallelized using mpi4py library, before being ported to C language by using Cython and compiled on TARA using foss-2021b toolchain.
The weak-scaling parallel efficiency $R(1)/R(N_\mathrm{s})$, where $R(N_\mathrm{s})$ is wall-clock runtime of the $N_\mathrm{s}$ scaled case, is plotted in Figure \ref{fig:WeakSca} from $N_\mathrm{s}=1$ to $N_\mathrm{s}=10$.
The results of the serial runtime $R(1)$ are provided in Table \ref{table:weaksca_base}.

From Figure \ref{fig:WeakSca}, the parallel efficiency of the integral-like schemes decreases as a larger time step size $\Delta t$ is chosen. 
However, from Table \ref{table:weaksca_base}, the wall-clock runtime for a larger time step size is evidently lower.
The parallel efficiency also declines as the problem size becomes larger and more CPU cores are used, but remains roughly unchanged at moderate-to-large scale cases. 
Therefore, the time step size of operational applications may need to be experimentally found to match available computational resources.

\begin{center}
\begin{table*}
	\centering
	\caption{The wall-clock serial runtime $R(1)$ of each integral-like scheme for different time step size $\Delta t$, but ran for the same simulation time, i.e., from $t_0 = 10$ to $t=810$.}
	\begin{tabular}{cccrrr}
		\hline
        \multirow{2}{*}{$\Delta\mathbf{t}$} & \multirow{2}{*}{$\mathbf{5\sigma}$} & \multirow{2}{*}{$\mathbf{d}$} & \multicolumn{3}{c}{Runtime (sec)} \\
        \cline{4-6}
		 &  & & LG & CG & QG \\
		\hline
		$0.2$ &  $3.16$ & $0.26-0.31$ & $1,945$ & $3,171$ & $4,714$ \\
		$0.8$ &  $6.32$ & $1.05-1.23$ &   $890$ & $1,420$ & $2,069$ \\
		$3.2$ & $12.65$ & $4.22-4.92$ &   $327$ &   $508$ &   $757$ \\
		\hline
	\end{tabular}
	\label{table:weaksca_base}
\end{table*}
\end{center}

\newpage

\textbf{Example 4.} Burgers'-Fisher equation
\begin{equation}
	\frac{\partial u}{\partial t} + u \frac{\partial u}{\partial x} = \nu \frac{\partial^2 u}{\partial x^2} - 3 u (1-u) (1-2u)
	\label{eq:burger_source}
\end{equation}
is considered to show a possible extension of the integral-like method.
Applying the split approach, which separates terms in the equation into stages and successively solve them, to Eq. (\ref{eq:burger_source}), the stage equations are
\begin{align}
	\frac{\partial u}{\partial t} + u \frac{\partial u}{\partial x} 
	&= \nu \frac{\partial^2 u}{\partial x^2} \label{eq:burger_source_split_1}\\
	\frac{\partial u}{\partial t} &= - 3 u (1-u) (1-2u)
\label{eq:burger_source_split_2}
\end{align}
The first stage, Eq. (\ref{eq:burger_source_split_1}), is to solve the Burger's equation; therefore, the procedures discussed in previous sections are directly employed.
The second stage, Eq. (\ref{eq:burger_source_split_2}), is a growth/decay equation.
Its analytical solution is found by solving
\begin{align*}
	-3 \, \mathrm{d}t = \Bigg(\, \frac{1}{u}-\frac{1}{(1-u)}+\frac{4}{(1-2u)} \, \Bigg)\, \mathrm{d}u
	 \quad \rightarrow \quad
	 A \, \mathrm{e}^{-3t} = \frac{1}{4} - \frac{1}{4 (2u-1)^2}
\end{align*}
Therefore, the integral-like scheme for this second stage is
\begin{equation*}
	u(x_i, t+\Delta t) = \frac{1}{2} \Big(1 \pm \frac{1}{\sqrt{1 - 4 A \, \mathrm{exp}(-3 \Delta t)}}\Big) 
\qquad \text{where} \quad
	A = \frac{1}{4} \bigg(1 - \frac{1}{(2u(x_i,t)-1)^2}\bigg)
\end{equation*}

\begin{figure}[h]
	\includegraphics[width=\linewidth]{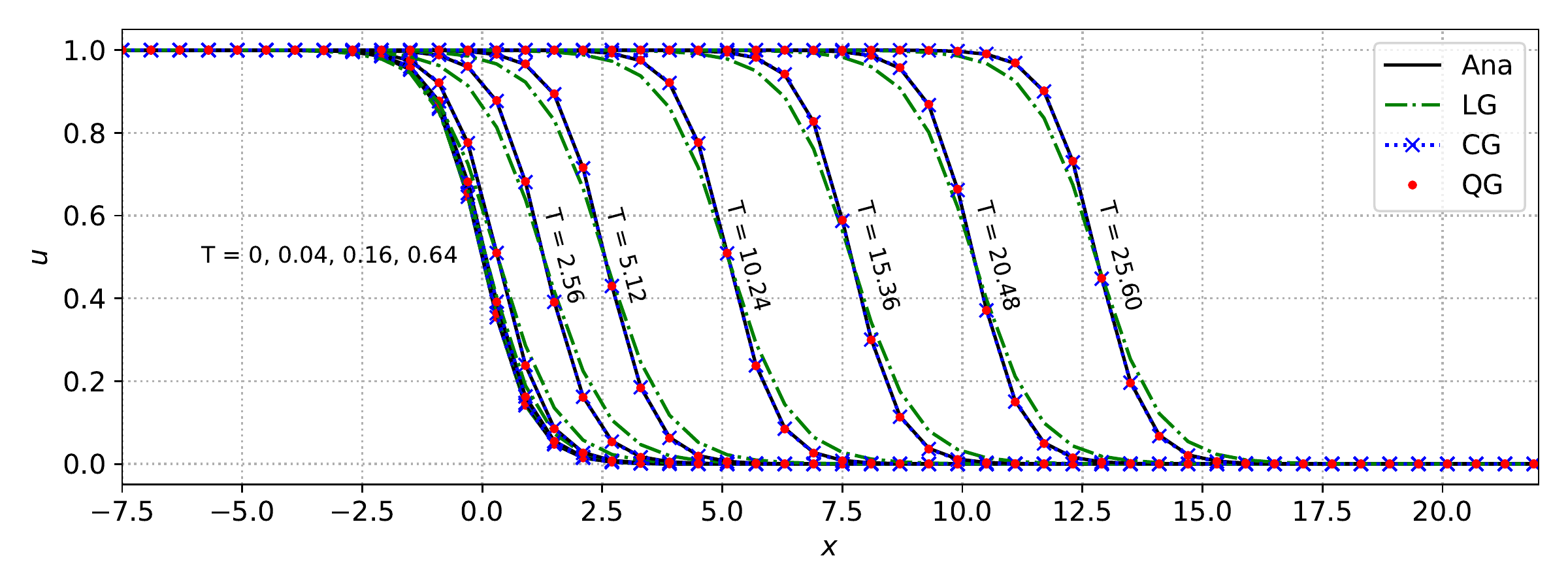}
	\caption{Numerical solutions and exact solution of Example 4: $\nu = 1$, $\Delta x = 0.6$ and $\Delta t = 0.04$}
	\label{fig:Sol_additionalsource_travel}
\end{figure}

From \cite{ramos2009}, an exact solution of Eq. (\ref{eq:burger_source}) when $\nu = 1$ is 
\begin{equation}
	u(x,t) = \frac{1}{2} \Big(1-\mathrm{tanh}\Big(x-\frac{t}{2}\Big)\Big)
	\label{eq:burger_source_sol}
\end{equation}
Using the same setup as in Figure \ref{fig:Sol_tophat_travel1} of Example 3, the numerical results of CG and QG agree well with the exact solution as seen in Figure \ref{fig:Sol_additionalsource_travel}.
The integral-like approach, therefore, has the potential for solving other Burger's-type equations as well.

\section{Conclusions}
\label{Sec:Conclusion}

In this study, a numerical approach based on the continuous representation of variables using local polynomial interpolation is explored.
When applied to solve the one-dimensional Burgers' equation, the schemes derived using linear (LG), cubic (CG) and quintic (QG) interpolations are found to be numerically stable at large time step size $\Delta t$, with a stability condition of $\nu \Delta t > 0.02 \Delta x^2$.
However, the order of accuracy is not consistently improved when a smaller grid spacing is employed; nonetheless, the smallest possible error norms are around the order of $10^{-4}$.
Being an explicit scheme, their weak-scaling parallel efficiency is generally adequate.
This approach shows promise for operational applications that favor reliability, fast computation and good parallel scalability, such as numerical weather prediction.

\section*{Acknowledgments}
I would like to thank NSTDA supercomputer center (ThaiSC) for providing computational resources for this work.

\bibliographystyle{elsarticle-num} 
\bibliography{IntegralLike}

\begin{thebibliography}{10}
\expandafter\ifx\csname url\endcsname\relax
  \def\url#1{\texttt{#1}}\fi
\expandafter\ifx\csname urlprefix\endcsname\relax\def\urlprefix{URL }\fi
\expandafter\ifx\csname href\endcsname\relax
  \def\href#1#2{#2} \def\path#1{#1}\fi

\bibitem{reviewBurger}
M.~P. Bonkile, A.~Awasthi, C.~Lakshmi, V.~Mukundan, V.~S. Aswin, A systematic
  literature review of burgers’ equation with recent advances, Pramana 90
  (2018) 1--21.
\newblock \href {https://doi.org/https://doi.org/10.1007/s12043-018-1559-4}
  {\path{doi:https://doi.org/10.1007/s12043-018-1559-4}}.

\bibitem{abdul2022}
M.~Abdullah, M.~Yaseen, M.~De~la Sen, {An efficient collocation method based on
  Hermite formula and cubic B-splines for numerical solution of the Burgers’
  equation}, Math. Comput. Simulation 197~(C) (2022) 166--184.
\newblock \href {https://doi.org/10.1016/j.matcom.2022.02.}
  {\path{doi:10.1016/j.matcom.2022.02.}}

\bibitem{dogan2004}
A.~Dogan, A galerkin finite element approach to burgers' equation, Appl. Math.
  Comput. 157~(2) (2004) 331--346.
\newblock \href {https://doi.org/https://doi.org/10.1016/j.amc.2003.08.037}
  {\path{doi:https://doi.org/10.1016/j.amc.2003.08.037}}.

\bibitem{gana2014}
I.~Ganaie, V.~Kukreja, Numerical solution of burgers’ equation by cubic
  hermite collocation method, Appl. Math. Comput. 237 (2014) 571--581.
\newblock \href {https://doi.org/https://doi.org/10.1016/j.amc.2014.03.102}
  {\path{doi:https://doi.org/10.1016/j.amc.2014.03.102}}.

\bibitem{hon1998}
Y.~Hon, X.~Mao, An efficient numerical scheme for burgers' equation, Appl.
  Math. Comput. 95~(1) (1998) 37--50.
\newblock \href {https://doi.org/https://doi.org/10.1016/S0096-3003(97)10060-1}
  {\path{doi:https://doi.org/10.1016/S0096-3003(97)10060-1}}.

\bibitem{huang2010}
P.~Huang, A.~Abduwali, The modified local crank–nicolson method for one- and
  two-dimensional burgers’ equations, Comput. Math. Appl. 59~(8) (2010)
  2452--2463.
\newblock \href {https://doi.org/https://doi.org/10.1016/j.camwa.2009.08.069}
  {\path{doi:https://doi.org/10.1016/j.camwa.2009.08.069}}.

\bibitem{aswin2017}
A.~Vs, A.~Awasthi, A differential quadrature based numerical method for highly
  accurate solutions of burgers' equation: Dqm based numerical method for
  burgers' equation, Numer. Meth. Part. D. E. 33 (07 2017).
\newblock \href {https://doi.org/10.1002/num.22178}
  {\path{doi:10.1002/num.22178}}.

\bibitem{yang2021}
X.~Yang, Y.~Ge, B.~Lan, A class of compact finite difference schemes for
  solving the 2d and 3d burgers’ equations, Math. Comput. Simulation 185
  (2021) 510--534.
\newblock \href {https://doi.org/https://doi.org/10.1016/j.matcom.2021.01.009}
  {\path{doi:https://doi.org/10.1016/j.matcom.2021.01.009}}.

\bibitem{guo2016}
Y.~Guo, Y.~feng Shi, Y.~min Li, A fifth-order finite volume weighted compact
  scheme for solving one-dimensional burgers’ equation, Appl. Math. Comput.
  281 (2016) 172--185.
\newblock \href {https://doi.org/https://doi.org/10.1016/j.amc.2016.01.061}
  {\path{doi:https://doi.org/10.1016/j.amc.2016.01.061}}.

\bibitem{gupta2021}
S.~Gupta, V.~K. Kukreja, An improvised collocation algorithm with specific end
  conditions for solving modified burgers equation, Numer. Meth. Part. D. E.
  37~(1) (2021) 874--896.
\newblock \href {https://doi.org/https://doi.org/10.1002/num.22557}
  {\path{doi:https://doi.org/10.1002/num.22557}}.

\bibitem{jena2023}
S.~R. Jena, G.~S. Gebremedhin, Decatic b-spline collocation scheme for
  approximate solution of burgers' equation, Numer. Meth. Part. D. E. 39~(3)
  (2023) 1851--1869.
\newblock \href {https://doi.org/https://doi.org/10.1002/num.22747}
  {\path{doi:https://doi.org/10.1002/num.22747}}.

\bibitem{jiang2021}
Y.~Jiang, X.~Chen, R.~Fan, X.~Zhang, High order semi-implicit weighted compact
  nonlinear scheme for viscous burgers’ equations, Math. Comput. Simulation
  190 (2021) 607--621.
\newblock \href {https://doi.org/https://doi.org/10.1016/j.matcom.2021.06.006}
  {\path{doi:https://doi.org/10.1016/j.matcom.2021.06.006}}.

\bibitem{mohan2015}
R.~K. Mohanty, J.~Talwar, A new compact alternating group explicit iteration
  method for the solution of nonlinear time-dependent viscous burgers’
  equation, Numer. Anal. Appl. 8 (2015) 314--328.
\newblock \href {https://doi.org/https://doi.org/10.1134/S1995423915040059}
  {\path{doi:https://doi.org/10.1134/S1995423915040059}}.

\bibitem{zhang2011}
R.~Zhang, Y.~Xi-Jun, Z.~Guo-Zhong, Local discontinuous galerkin method for
  solving burgers and coupled burgers equations, Chin. Phys. B 20~(11) (2011)
  110205.
\newblock \href {https://doi.org/10.1088/1674-1056/20/11/110205}
  {\path{doi:10.1088/1674-1056/20/11/110205}}.

\bibitem{kadal2006}
M.~K. Kadalbajoo, A.~Awasthi, A numerical method based on crank-nicolson scheme
  for burgers’ equation, Appl. Math. Comput. 182~(2) (2006) 1430--1442.
\newblock \href {https://doi.org/https://doi.org/10.1016/j.amc.2006.05.030}
  {\path{doi:https://doi.org/10.1016/j.amc.2006.05.030}}.

\bibitem{kan2012}
R.~Kannan, Z.~Wang, A high order spectral volume solution to the burgers'
  equation using the hopf–cole transformation, Internat. J. Numer. Methods
  Fluids 69~(4) (2012) 781--801.
\newblock \href {https://doi.org/https://doi.org/10.1002/fld.2612}
  {\path{doi:https://doi.org/10.1002/fld.2612}}.

\bibitem{kum2019}
S.~S. Kumbhar, S.~Thakar, Galerkin finite element method for forced burgers'
  equation, J. Math. Model. 7~(4) (2019) 445--467.
\newblock \href {https://doi.org/10.22124/jmm.2019.13259.1265}
  {\path{doi:10.22124/jmm.2019.13259.1265}}.

\bibitem{kut1999}
S.~Kutluay, A.~Bahadir, A.~Özdeş, Numerical solution of one-dimensional
  burgers equation: explicit and exact-explicit finite difference methods, J.
  Comput. Appl. Math. 103~(2) (1999) 251--261.
\newblock \href {https://doi.org/https://doi.org/10.1016/S0377-0427(98)00261-1}
  {\path{doi:https://doi.org/10.1016/S0377-0427(98)00261-1}}.

\bibitem{liao2008}
W.~Liao, An implicit fourth-order compact finite difference scheme for
  one-dimensional burgers’ equation, Appl. Math. Comput. 206~(2) (2008)
  755--764.
\newblock \href {https://doi.org/https://doi.org/10.1016/j.amc.2008.09.037}
  {\path{doi:https://doi.org/10.1016/j.amc.2008.09.037}}.

\bibitem{mukun2015}
V.~Mukundan, A.~Awasthi, Efficient numerical techniques for burgers’
  equation, Appl. Math. Comput. 262 (2015) 282--297.
\newblock \href {https://doi.org/https://doi.org/10.1016/j.amc.2015.03.122}
  {\path{doi:https://doi.org/10.1016/j.amc.2015.03.122}}.

\bibitem{pandy2009}
K.~Pandey, L.~Verma, A.~K. Verma, On a finite difference scheme for burgers’
  equation, Appl. Math. Comput. 215~(6) (2009) 2206--2214.
\newblock \href {https://doi.org/https://doi.org/10.1016/j.amc.2009.08.018}
  {\path{doi:https://doi.org/10.1016/j.amc.2009.08.018}}.

\bibitem{sakai2005}
K.~Sakai, I.~Kimura, A numerical scheme based on a solution of nonlinear
  advection–diffusion equations, J. Comput. Appl. Math. 173~(1) (2005)
  39--55.
\newblock \href {https://doi.org/https://doi.org/10.1016/j.cam.2004.02.019}
  {\path{doi:https://doi.org/10.1016/j.cam.2004.02.019}}.

\bibitem{xie2008}
S.-S. Xie, S.~Heo, S.~Kim, G.~Woo, S.~Yi, Numerical solution of one-dimensional
  burgers’ equation using reproducing kernel function, J. Comput. Appl. Math.
  214~(2) (2008) 417--434.
\newblock \href {https://doi.org/https://doi.org/10.1016/j.cam.2007.03.010}
  {\path{doi:https://doi.org/10.1016/j.cam.2007.03.010}}.

\bibitem{zhao2011}
G.~Zhao, X.~Yu, R.~Zhang, The new numerical method for solving the system of
  two-dimensional burgers’ equations, Comput. Math. Appl. 62~(8) (2011)
  3279--3291.
\newblock \href {https://doi.org/https://doi.org/10.1016/j.camwa.2011.08.044}
  {\path{doi:https://doi.org/10.1016/j.camwa.2011.08.044}}.

\bibitem{aksan2005}
E.~Aksan, A numerical solution of burgers’ equation by finite element method
  constructed on the method of discretization in time, Appl. Math. Comput.
  170~(2) (2005) 895--904.
\newblock \href {https://doi.org/https://doi.org/10.1016/j.amc.2004.12.027}
  {\path{doi:https://doi.org/10.1016/j.amc.2004.12.027}}.

\bibitem{cal1982}
J.~Caldwell, P.~Smith, Solution of burgers' equation with a large reynolds
  number, Appl. Math. Model. 6~(5) (1982) 381--385.
\newblock \href {https://doi.org/https://doi.org/10.1016/S0307-904X(82)80102-9}
  {\path{doi:https://doi.org/10.1016/S0307-904X(82)80102-9}}.

\bibitem{cal1981}
J.~Caldwell, P.~Wanless, A.~Cook, A finite element approach to burgers'
  equation, Appl. Math. Model. 5~(3) (1981) 189--193.
\newblock \href {https://doi.org/https://doi.org/10.1016/0307-904X(81)90043-3}
  {\path{doi:https://doi.org/10.1016/0307-904X(81)90043-3}}.

\bibitem{chai2020}
Y.~Chai, J.~Ouyang, Appropriate stabilized galerkin approaches for solving
  two-dimensional coupled burgers’ equations at high reynolds numbers,
  Comput. Math. Appl. 79~(5) (2020) 1287--1301.
\newblock \href {https://doi.org/https://doi.org/10.1016/j.camwa.2019.08.036}
  {\path{doi:https://doi.org/10.1016/j.camwa.2019.08.036}}.

\bibitem{arora2013}
G.~Arora, B.~K. Singh, Numerical solution of burgers’ equation with modified
  cubic b-spline differential quadrature method, Appl. Math. Comput. 224 (2013)
  166--177.
\newblock \href {https://doi.org/https://doi.org/10.1016/j.amc.2013.08.071}
  {\path{doi:https://doi.org/10.1016/j.amc.2013.08.071}}.

\bibitem{ghase2018}
M.~Ghasemi, An efficient algorithm based on extrapolation for the solution of
  nonlinear parabolic equations, Int. J. Nonlinear Sci. Numer. Simul. 19~(1)
  (2018) 37--51.
\newblock \href {https://doi.org/doi:10.1515/ijnsns-2017-0060}
  {\path{doi:doi:10.1515/ijnsns-2017-0060}}.

\bibitem{singh2021}
B.~K. Singh, M.~Gupta, A new efficient fourth order collocation scheme for
  solving burgers’ equation, Appl. Math. Comput. 399 (2021) 126011.
\newblock \href {https://doi.org/https://doi.org/10.1016/j.amc.2021.126011}
  {\path{doi:https://doi.org/10.1016/j.amc.2021.126011}}.

\bibitem{tam2016}
M.~Tamsir, N.~Dhiman, V.~K. Srivastava, Extended modified cubic b-spline
  algorithm for nonlinear burgers' equation, Beni-Suef Univ. J. Basic Appl.
  Sci. 5~(3) (2016) 244--254.
\newblock \href {https://doi.org/https://doi.org/10.1016/j.bjbas.2016.09.001}
  {\path{doi:https://doi.org/10.1016/j.bjbas.2016.09.001}}.

\bibitem{gao2013}
Y.~Gao, L.-H. Le, B.-C. Shi, Numerical solution of burgers’ equation by
  lattice boltzmann method, Appl. Math. Comput. 219~(14) (2013) 7685--7692.
\newblock \href {https://doi.org/https://doi.org/10.1016/j.amc.2013.01.056}
  {\path{doi:https://doi.org/10.1016/j.amc.2013.01.056}}.

\bibitem{kumar2020}
N.~Kumar, R.~Majumdar, S.~Singh, Predictor–corrector nodal integral method
  for simulation of high reynolds number fluid flow using larger time steps in
  burgers ’ equation, Comput. Math. Appl. 79~(5) (2020) 1362--1381.
\newblock \href {https://doi.org/https://doi.org/10.1016/j.camwa.2019.09.001}
  {\path{doi:https://doi.org/10.1016/j.camwa.2019.09.001}}.

\bibitem{lara2022}
F.~M. {de Lara}, E.~Ferrer, Accelerating high order discontinuous galerkin
  solvers using neural networks: 1d burgers’ equation, Comput. \& Fluids 235
  (2022) 105274.
\newblock \href
  {https://doi.org/https://doi.org/10.1016/j.compfluid.2021.105274}
  {\path{doi:https://doi.org/10.1016/j.compfluid.2021.105274}}.

\bibitem{fluidgraphic}
R.~Bridson, Fluid Simulation for Computer Graphics, Second Edition, Taylor \&
  Francis, 2015.

\bibitem{olver2013}
P.~Olver, Introduction to Partial Differential Equations, Undergraduate Texts
  in Mathematics, Springer International Publishing, 2013.

\bibitem{stone2009}
M.~Stone, P.~Goldbart, Mathematics for Physics: A Guided Tour for Graduate
  Students, Cambridge University Press, 2009.

\bibitem{Sarbo2014}
M.~Sarboland, A.~Aminataei, On the numerical solution of one-dimensional
  nonlinear nonhomogeneous burgers' equation, J. Appl. Math. 2014 (2014)
  598432:1--598432:15.

\bibitem{benton1972}
E.~R. Bentom, G.~W. Platzman, A table of solutions of the one-dimensional
  burgers equation, Quart. Appl. Math. 30~(2) (1972) 195--212.

\bibitem{ramos2009}
J.~Ramos, Picard’s iterative method for nonlinear
  advection–reaction–diffusion equations, Appl. Math. Comput. 215~(4)
  (2009) 1526--1536.
\newblock \href {https://doi.org/https://doi.org/10.1016/j.amc.2009.07.004}
  {\path{doi:https://doi.org/10.1016/j.amc.2009.07.004}}.

\end{thebibliography}


\end{document}